\documentclass[12pt]{extarticle}
\IfFileExists{meta/preamble.sty}{\usepackage{meta/preamble}}{%
  \usepackage[T1]{fontenc}
  \usepackage{lmodern,amsmath,amssymb,amsthm,mathtools}
  \usepackage{tikz-cd}
  \setlength{\emergencystretch}{1.5em}
  \usepackage[margin=1.1in]{geometry}
  \usepackage[hidelinks]{hyperref}
  \newtheorem{theorem}{Theorem}[section]
  \newtheorem{proposition}[theorem]{Proposition}
  \newtheorem{lemma}[theorem]{Lemma}
  \newtheorem{corollary}[theorem]{Corollary}
  \theoremstyle{definition}
  \newtheorem{definition}[theorem]{Definition}
  \newtheorem{construction}[theorem]{Construction}
  \newtheorem{discussion}[theorem]{Discussion}
  \newtheorem{remark}[theorem]{Remark}
}
\usepackage{hyperref}
\hypersetup{pdftitle={E-infinity-Coherent Motivic Chern Classes},pdfauthor={Amos Kaminski}}
\providecommand{\SH}{\mathcal{SH}}
\providecommand{\MGL}{\mathrm{MGL}}
\providecommand{\Pic}{\operatorname{Pic}}
\providecommand{\Fun}{\operatorname{Fun}}
\providecommand{\CAlg}{\operatorname{CAlg}}
\providecommand{\Mod}{\operatorname{Mod}}
\providecommand{\Map}{\operatorname{Map}}

\providecommand{\B}{\mathcal{B}}
\providecommand{\C}{\mathcal{C}}

\theoremstyle{plain}
\newtheorem{introtheorem}{Theorem}

\newtheorem{introdefinition}{Definition}

\title{ \(\mathbb{E}_\infty\)-Coherent Motivic Chern Classes}
\author{Amos Kaminski}
\date{}

\begin{document}
\maketitle

\begin{abstract}
In \cite[Corollary~4.5]{SL} Carmeli and Luecke constructed a
Tate-valued $\mathbb E_\infty$-lift of the total Chern class
of complex vector bundles. In this paper, we define a Tate
construction for genuine $\mathbb G_m$-equivariant motivic spectra
and use it to construct a motivic analogue of their lift. For
integral motivic cohomology, we show that this lift recovers the
motivic Chern classes.
\end{abstract}

\par
\vspace*{\fill}

\noindent
\begin{minipage}{\textwidth}
  \centering
  \includegraphics[width=\linewidth]{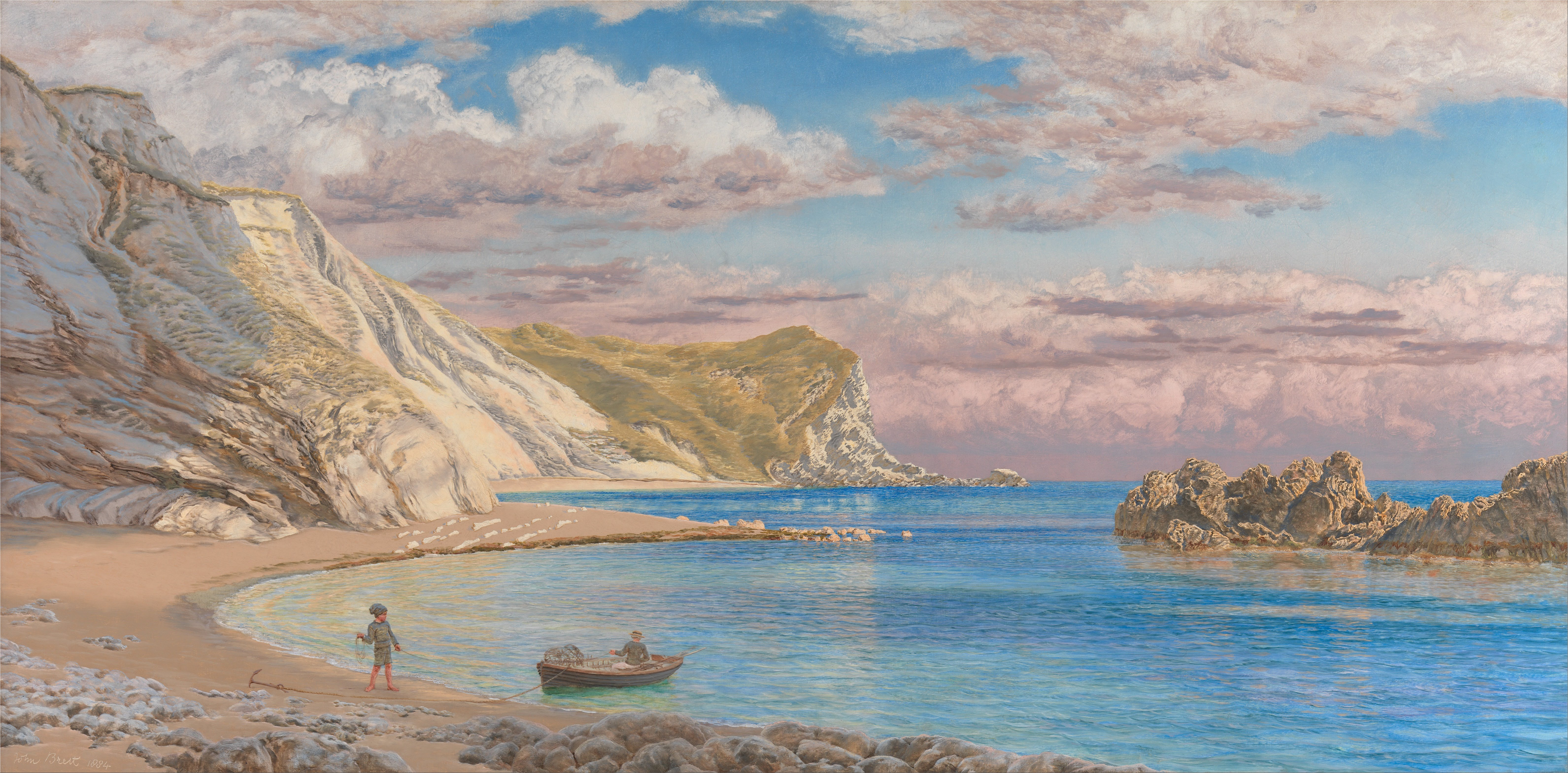}
  \par\smallskip
  {\footnotesize John Brett, \emph{Man of War Rocks, Coast of Dorset} (1884).\par}
\end{minipage}
\par
\vspace*{\fill}
\clearpage
\newpage
\tableofcontents
\newpage
\section*{Introduction}

Chern classes are fundamental invariants in both algebraic topology and algebraic geometry.
A key property shared by the two settings is their \emph{multiplicativity}.

\begin{theorem}[Whitney sum formula, {\cite[\S14, formula~(14.7)]{MS}, \cite[\S8, before Lemma~8.3]{AHI}}]
Let \(V\) and \(V'\) be complex vector bundles over a paracompact Hausdorff space or algebraic vector bundles over a qcqs scheme \(S\), and denote by \(c(V)\) the total Chern class. Then
\[
c(V \oplus V') = c(V)\,c(V'),
\]
where the multiplication takes place in the ordinary or motivic cohomology ring of the base.
\end{theorem}

In the topological case, direct sum gives $BU$ its $H$-space structure.
The Whitney sum formula says that the total Chern class takes
direct sums to products. It therefore gives a map of $H$-spaces
\[
BU \longrightarrow \mathrm{GL}_1(H\mathbb Z[z]),
\]
where $z$ has homotopy degree $2$. This map sends a vector bundle to its total Chern class with parameter $z$ (see \cite[Example~2.11]{SL}).

A natural question then arises: can this construction be lifted coherently to the \(\mathbb{E}_\infty\) level?
This question goes back to Segal \cite{Segal}. Boyer, Lawson, Lima-Filho, Mann, and Michelsohn constructed an infinite loop lift using spaces of algebraic cycles \cite{BLLMM}.
Carmeli and Luecke took a different approach to obtaining an \(\mathbb{E}_\infty\)-coherent lift \cite[Definition~4.1 and Corollary~4.5]{SL}.

The first step is to understand how these two $H$-spaces can be lifted.
A natural lift for $BU$ is $\mathrm{bu}$. For $\mathrm{GL}_1(H\mathbb Z[z])$, there are several candidates:\\ The most natural one is the spectrum of units $\mathrm{gl}_1(H\mathbb Z[z])$ of the usual graded polynomial ring. Snaith showed that the total Chern class does not lift to a map of spectra with this target \cite{Snaith1977}. 

Hence Carmeli and Luecke choose another candidate, namely the spectrum of units of the Tate spectrum $(H\mathbb Z)^{t\mathbb T}$, where $\mathbb T=S^1$.

Carmeli and Luecke introduced a notion of \emph{projective orientation}
\cite[Definition~3.4]{SL}. For an $\mathbb E_\infty$-oriented ring
spectrum $(R,\omega)$, the Tate construction admits two natural
projective orientations. One is induced by $\omega$, and the other
is the Euler orientation.

\begin{definition}[Euler--Tate characteristic class, {\cite[Definition~4.1]{SL}}]
The \emph{Euler--Tate characteristic class} of $(R,\omega)$ is the
map of spectra obtained by dividing the orientation induced by $\omega$
by the Euler orientation:
\[
    c^{\mathrm{et}}\colon \mathrm{bu}\longrightarrow
    \mathrm{gl}_1(R^{t\mathbb T}).
\]
\end{definition}

Their main result is that, for ordinary integral cohomology, the Euler--Tate characteristic
class recovers the Chern classes. This gives another answer to Segal's question \cite{Segal}.

\begin{theorem}[{\cite[Corollary~4.5]{SL}}]
Let $R=H\mathbb Z$ with its canonical orientation, and let $t=c_1(L)$
for the tautological line bundle $L$ over $B\mathbb T$. For a complex
vector bundle $V$ of rank $r$, the Euler--Tate characteristic class is
\[
    c^{\mathrm{et}}([V]-[\mathbb C^r])=\sum_{i=0}^{r}c_i(V)\,t^{-i},
\]
where $\mathbb C^r$ denotes the trivial bundle of rank $r$.
\end{theorem}

Carmeli and Luecke ask whether the total motivic Chern class admits
an analogous coherent lift \cite[Question~6.2]{SL}. In this paper
, we answer this question in the affirmative by constructing such a lift,
with $\mathbb G_m$ playing the role of the circle:
The motivic analog uses the following motivic analog for the Tate construction:
\begin{introdefinition}[Motivic Tate construction]
For a genuine $\mathbb G_m$-equivariant motivic spectrum $D$ over a
scheme $S$ with trivial action, define its Tate construction by:
\[
    D^{t\mathbb G_m}
    := q_*\bigl(\widehat D[e_{\mathcal L}^{-1}]\bigr),
    \qquad q\colon S\times B\mathbb G_m\longrightarrow S.
\]
Here $\widehat {(-)} $ denotes Borel completion and $e_{\mathcal L}$ is
the Euler map of the tautological line bundle $\mathcal L$.
\end{introdefinition}

We define a motivic analogue for projective orientations,  and again for an $\mathbb E_\infty$-oriented ring
spectrum $(R,\omega)$, the Tate construction admits two natural
projective orientations: $\omega_{\mathrm{prj}},\  \omega_{\mathrm{eul}}$

\begin{introdefinition}[Motivic Euler--Tate characteristic class]\label{intro:euler-tate}
The \emph{motivic Euler--Tate characteristic class} of $(R,\omega)$
is the map of grouplike motivic $\mathbb E_\infty$-spaces
\[
    c^{\mathrm{et}}_{R,\omega}
    :=\frac{\omega_{\mathrm{prj}}}{\omega_{\mathrm{eul}}}
    \colon \mathrm{BGL}\longrightarrow\mathrm{GL}_1(R^{t\mathbb G_m}).
\]
\end{introdefinition}

And for integral motivic cohomology, our construction recovers the Chern
polynomial.

\begin{introtheorem}[Integral motivic Chern classes]\label{intro:integral}
Let $H\mathbb Z_S$ have its canonical orientation and write
$t=c_1^{H\mathbb Z}(\mathcal L)$. For a vector bundle $V$ of rank $r$
on a smooth $S$-scheme $X$,
\[
    c^{\mathrm{et}}_{H\mathbb Z}\bigl([V]-[\mathcal O_X^r]\bigr)
    =\sum_{i=0}^r c_i(V)\,t^{-i},
\]
where $c_i(V)$ is the $i$th motivic Chern class of $V$.
\end{introtheorem}

This polynomial gives a class in the Tate cohomology of $X$.
More generally, we compute the $R^{t\mathbb G_m}$-cohomology of $X$
in every bidegree:
\begin{introtheorem}[Tate coefficients]\label{intro:coefficients}
For each smooth qcqs $X$ over $S$ and each bidegree $(a,b)$, 
\[
 (R^{t\mathbb G_m})^{a,b}(X)
 \cong \underset{m\geq 0}{\operatorname{colim}}\;
       \prod_{j\geq -m} R^{a-2j,b-j}(X)\,t^j .
\]
where $t=c_1^R(\mathcal L)$.
\end{introtheorem}
\begin{remark}
Carmeli and Luecke also give an $\mathbb E_\infty$-refinement of the
sharp construction $\omega\mapsto\omega^{\#}$ of Ando, French, and
Ganter. Applying it to the complex string orientation gives an
$\mathbb E_\infty$-lift of the Jacobi orientation
\cite[Theorem~5.14 and Corollary~5.18]{SL}. We expect these results
to have motivic analogues and leave this question for future work.
\end{remark}
\subsection*{Acknowledgments}
First and foremost, I would like to express my heartfelt gratitude to my PhD advisor, Shachar Carmeli, for his guidance and support.
I am grateful to Yuval Lotenberg for helpful discussions on the Tate construction.
I would also like to thank Kiran Luecke and Akhil Mathew for the stimulating conversations we had during my visit to the University of Chicago.
I am also grateful to Noam Zimhoni and Guy Kapon for our conversations about internal categories.
Finally, I would like to thank everyone in my research group.

\subsection*{Funding}
This work was partially supported by BSF Grant No. 2024766 and by the Azrieli Foundation.
\subsection*{Acknowledgment of AI assistance}
The mathematical ideas and research direction are due to the authors, who wrote nearly all of the manuscript about a year ago. During the final revision, the authors used OpenAI’s ChatGPT and Codex to help fill in a few remaining technical details in the proofs, to improve phrasing, and to proofread the manuscript. The authors take full responsibility for the mathematical content.
\section{Motivic Spectra: Background and Notation}

We first recall the definitions of the unstable and stable category of motives over a base scheme $S$.

\subsection{Basic Definitions}

In motivic homotopy theory, we want the affine line $\mathbb A^1$ to plays
the role of the interval $[0,1]$. We therefore require a motivic
space $F$ to carry each projection $X\times\mathbb A^1\to X$ to an
equivalence. We also require $F$ to satisfy Nisnevich descent.

Let $S$ be a quasi-compact quasi-separated (qcqs) scheme. Write $\mathrm{Sm}_S$ for the site of smooth finitely presented $S$-schemes. 
\begin{definition}[Motivic Spaces, {\cite[\S3]{MV}, \cite[\S2.2]{norms}}]
The \emph{category of motivic spaces} over $S$ is defined as the intersection of Nisnevich sheaves and $\mathbb{A}^1$-invariant presheaves of spaces:
\[
\mathcal{H}(S) = \mathrm{Shv}_{\mathrm{Nis}}(\mathrm{Sm}_S) \cap \mathcal{P}_{\mathbb{A}^1}(\mathrm{Sm}_S) \subset \mathcal{P}(\mathrm{Sm}_S).
\]
\end{definition}

We now want to define motivic spectra by stabilizing the unstable motivic homotopy category. In this category, there are two natural candidates for the sphere. The first is the constant circle \(\underline{S^1}\), and the second is the multiplicative group \(\mathbb{G}_m\) pointed at \(1\). We may therefore suspend with respect to either of these spheres. The two indices of a motivic sphere keep track of these two kinds of suspension:

\begin{definition}[Motivic Spheres]
For $p\ge q\ge0$, the \emph{motivic sphere} is defined as
\[
S^{p,q} := \underline{S}^{p-q} \wedge \mathbb{G}_m^{\wedge q},
\]
where $\underline{S}^{p-q}$ is the constant $(p-q)$-sphere and $\mathbb{G}_m = \mathbb{A}^1 \setminus \{0\}$ is pointed at $1$.
\end{definition}

Stabilizing with respect to 
\[
S^{2,1}=\underline{S^1}\wedge(\mathbb G_m,1)\simeq(\mathbb P^1,\infty).
\]
makes suspension by both $S^1$ and $\mathbb G_m$ invertible, and
hence makes suspension by every motivic sphere invertible.
We will now define the category of motivic spectra:

\begin{definition}[Motivic Spectra, {\cite[\S4.1]{norms}}]
The \emph{category of motivic spectra} $\mathcal{SH}(S)$ over $S$
is obtained by inverting the motivic sphere $S^{2,1}$ in the
category of pointed motivic spaces. This yields a stable
$\infty$-category with a symmetric monoidal structure induced by the smash product.
\end{definition}

In $\mathcal{SH}(S)$ we use the notation $S^{p,q}$ for all
$p,q\in\mathbb Z$ and for a motivic spectrum $R$ over $S$ and $X\in\mathrm{Sm}_S$, write $R_X$ for its pullback to $X$ and
\[
 R^{a,b}(X)=[\mathbf1_X,\Sigma^{a,b}R_X]_{\SH(X)}.
\].

The invertible objects in $\mathcal{SH}(S)$ with respect to the monoidal structure on $\mathcal{SH}(S)$ form the Picard space of motivic spectra:

\begin{definition}[Sph, {\cite[\S16.1]{norms}}]
Let $\mathrm{Sph}(S) = \Pic(\mathcal{SH}(S))$ be the grouplike $\mathbb E_\infty$-space of tensor-invertible motivic spectra and their equivalences.
\end{definition}

\subsection{Motivic Colimits}
Let $\B$ be an $\infty$-topos. In the language of Martini and Wolf
\cite{MWcolim}, a category internal to $\B$ can be described as a
sheaf of $\infty$-categories on $\B$, that is, a functor
$\B^{\mathrm{op}}\to\mathrm{Cat}_\infty$ that sends colimits in
$\B$ to limits. An object $I\in\B$ determines a global object of
the universe $\Omega_{\B}$ of small internal groupoids. We may
therefore use $I$ to index diagrams in a $\B$-category and take
their internal colimits whenever they exist.

We now specialize to $\B=\mathcal P(\mathrm{Sm}_S)$. In this case,
a $\B$-category is determined by a presheaf of $\infty$-categories
on $\mathrm{Sm}_S$, and an object $I\in\B$ is a presheaf of spaces
on $\mathrm{Sm}_S$. In particular, $X\mapsto\SH(X)$ defines a
$\B$-category. The internal colimit of an $I$-indexed diagram in
$\SH(-)$ is its motivic colimit:

\begin{construction}[Motivic colimits, {\cite[Construction~2.3 and Definition~2.4]{mcolimits}}]
Let $\C$ be a $\mathcal P(\mathrm{Sm}_S)$-category admitting the indicated internal
colimits, and write $\pi:I\to 1$ for the projection. The
\emph{motivic colimit} of a diagram $f:I\to\C$ is the internal
colimit $\pi_!f\in\C(S)$. Equivalently, it is characterized by
\[
    \Map_{\C(S)}(\pi_!f,E)
    \simeq \Map_{\Fun_{\B}(I,\C)(1)}(f,\pi^*E).
\]
Concretely, for $\C=\SH(-)$ this gives
\[
    \pi_!f\simeq
    \operatorname*{colim}_{(p:X\to S,\,i\in I(X))}
       p_\# f_X(i).
\]
The colimit on the right is taken in $\SH(S)$ over the category
of elements of $I$. A morphism from $(X,i)$ to $(Y,j)$ consists of
a map $g:X\to Y$ over $S$ and an identification $i\simeq g^*j$.
In particular, for a family represented by a smooth $S$-scheme
$p:X\to S$ with value $E\in\SH(X)$, the motivic colimit is $p_\#E$.
\end{construction}

When the objects of the diagram are invertible, its motivic colimit
is the Thom spectrum of the family.
\begin{definition}[Motivic Thom spectrum, {\cite[Definition~16.1]{norms}}]
Let $I$ be a presheaf of spaces on $\mathrm{Sm}_S$. And let $f\colon I\to\mathrm{Sph}$ be a family of
invertible motivic spectra parametrized by $I$. Its \emph{motivic Thom spectrum}
$M(f)\in\SH(S)$ is the motivic colimit of the composite
\[
    I\xrightarrow{f}\mathrm{Sph}\hookrightarrow\SH(-).
\]
\end{definition}
\subsection{Motivic Cobordism and Orientations}
We now define the motivic analogue of the $J$-homomorphism.
Recall that $K$ denotes the presheaf on $\mathrm{Sm}_S$ that assigns
to each $X$ the algebraic $K$-theory space of perfect complexes on
$X$. Direct sum gives $K(X)$ its grouplike $\mathbb E_\infty$-structure
\cite[\S16.2]{norms}.
\begin{construction}[Motivic $J$-homomorphism, {\cite[\S16.2]{norms}}]
For a vector bundle $V$ over $X$, let
\[
    S^V:=\Sigma^\infty_X\bigl(V/(V\setminus 0_X)\bigr).
\]
The assignment $V\mapsto S^V$ is compatible with pullback and
satisfies $S^{V\oplus W}\simeq S^V\otimes S^W$. It extends to the
\emph{motivic $J$-homomorphism}
\[
    J:K\longrightarrow\mathrm{Sph}.
\]
For a virtual bundle $[V]-[W]$ over $X$, its value is
\[
    J_X([V]-[W])\simeq S^V\otimes (S^W)^{-1}.
\]
\end{construction}

Restricting $J$ to virtual bundles of rank zero gives the family
whose Thom spectrum is the motivic cobordism spectrum $MGL$:

\begin{construction}[MGL, {\cite[Theorem~16.13]{norms}}]
Let $e:K^\circ\to K$ denote the inclusion of rank $0$ vector bundles. The \emph{motivic algebraic cobordism spectrum} is the Thom spectrum
\[
    \MGL_S\simeq M(J|_{K^\circ}).
\]
\end{construction}

\section{Orientation of \texorpdfstring{$E_\infty$}{E-infinity} Motivic Spectra}

Let $R$ be an $\mathbb E_\infty$-ring spectrum. An
\emph{$\mathbb E_\infty$-orientation} of $R$ is a map of
$\mathbb E_\infty$-ring spectra
\[
    MU \longrightarrow R.
\]
In stable motivic homotopy theory, algebraic cobordism $\MGL_S$
plays the role of complex cobordism $MU$. This leads to the
following definition:

\begin{definition}
Let $R$ be a motivic $\mathbb E_\infty$-ring spectrum in $\SH(S)$.
An \emph{$\mathbb E_\infty$-orientation} of $R$ is a map of
$\mathbb E_\infty$-ring spectra
\[
\MGL_S \longrightarrow R.
\]
\end{definition}

\subsection{Characterization via Null Homotopies}

For ordinary Thom spectra, the universal property of
\cite[Theorem~3.5 and Lemma~3.15]{ACB} identifies a map of
$\mathbb E_\infty$-ring spectra $M(f)\to A$ with a multiplicative
trivialization of the family of invertible $A$-modules
$x\mapsto A\otimes f(x)$. We prove the motivic analogue below.
Applied to $\MGL_S$, it identifies an $\mathbb E_\infty$-orientation
of $R$ with a choice of equivalences
\[
    R_U\otimes S^\xi \simeq R_U
\]
for rank-zero virtual bundles $\xi$ on smooth $S$-schemes $U$,
natural in $\xi$, compatible with pullback in $U$, and symmetric
monoidal under direct sum.

Let $I$ be a presheaf of grouplike $\mathbb E_\infty$-spaces on
$\mathrm{Sm}_S$, and let $f:I\to\mathrm{Sph}$ be a map of such
presheaves, with motivic Thom spectrum $M(f)$. For
$A\in\CAlg(\SH(S))$, let
\[
 f_A:I\xrightarrow{f}\mathrm{Sph}\longrightarrow\Pic(\Mod_A)
\]
be its extension of scalars.

\begin{proposition}
\label{prop:orientation}
The space of maps of $\mathbb E_\infty$-ring spectra $M(f)\to A$
is equivalent to the space of null homotopies of $f_A$ as a map of
presheaves of grouplike $\mathbb E_\infty$-spaces.
\end{proposition}
We first explain the general idea behind the proof, leaving the details to
Section~\ref{proofof}.
By the universal property of the internal colimit defining $M(f)$,
a map of motivic spectra $M(f)\to A$ amounts to a family of maps
$f(x)\to A_U$, natural in $x\in I(U)$ and compatible with pullback
in $U$. To prove Proposition~\ref{prop:orientation}, we carry out
the following argument coherently. A map $f(x)\to A_U$ is
equivalently a map of $A_U$-modules $A_U\otimes f(x)\to A_U$.
For a multiplicative family, this map is an equivalence, since
tensoring it over $A_U$ with the corresponding map for $-x$ gives
the map for $x+(-x)=0$, which is the identity of $A_U$. For transforming this intuition into a proof we will need the following lemmas.

\subsection{Preparatory lemmas}

Equip $\mathcal P(\mathrm{Sm}_S)$ with its cartesian symmetric
monoidal structure. The presheaf
\[
 X\longmapsto \mathrm{Sph}(X)=\Pic(\SH(X))
\]
is a commutative algebra object in this category hence by \cite[Remark~2.2.2.5]{LurieHA} $\mathcal P(\mathrm{Sm}_S)/_\mathrm{Sph}$ inherits a symmetric monoidal structure:
Explicitly, for a finite set $J$ and families
$f_j:I_j\to\mathrm{Sph}$, their tensor product is the family
\[
 \prod_{j\in J} I_j
 \xrightarrow{\prod_{j\in J}f_j}
 \mathrm{Sph}^{J}
 \xrightarrow{\mu_J}
 \mathrm{Sph},
\]
where $\mu_J$ is the $J$-ary multiplication of $\mathrm{Sph}$.
For $J=\varnothing$, this gives the unit
$yS\simeq *\to\mathrm{Sph}$, whose value at $X$ selects
$\mathbf1_X\in\Pic(\SH(X))$.

\begin{lemma}[Symmetric monoidality of the Thom functor]
\label{lem:thom-families-monoidal}
The motivic Thom functor
\[
 M:\mathcal P(\mathrm{Sm}_S)_{/\mathrm{Sph}}\longrightarrow\SH(S)
\]
is symmetric monoidal.
\end{lemma}

\begin{proof}
Let $\int_{\mathrm{Sm}_S}\SH(-)$ denote the Grothendieck construction
of $X\mapsto\SH(X)$ with transition functors given by pullback.
Its objects are pairs $(X,E)$ with $X\in\mathrm{Sm}_S$ and $E\in\SH(X)$.
A morphism $(X,E)\to(Y,F)$ consists of an $S$-morphism $h:X\to Y$
and a map $E\to h^*F$ in $\SH(X)$.
By \cite[Theorem~4.5(ii)--(iii)]{KR} $\SH(-): \mathrm{Sm}_S^{\mathrm{op}}\to\CAlg(\mathrm{Cat}_\infty)$.
By Drew and Gallauer \cite[Corollary~A.12 and Remark~A.13(i)]{DrewGallauer},
the total category carries a symmetric monoidal structure\footnote{For our
convention on morphisms, take the opposite of the total category described
in their Remark~2.10.} with external tensor product
\[
(X,E)\boxtimes(Y,F)
=
\bigl(X\times_S Y,\,
\mathrm{pr}_X^*E\otimes\mathrm{pr}_Y^*F\bigr),
\]
with unit $(S,\mathbf1_S)$.
In every arity the tensor is given by
pullback along the projections from the corresponding fiber product,
followed by tensor product there.

By functoriality of the symmetric monoidal Grothendieck construction,
the inclusion of the fiber over the terminal object $S$ is symmetric monoidal:
\[
i:\SH(S)\hookrightarrow\int_{\mathrm{Sm}_S}\SH(-),
\qquad B\longmapsto(S,B).
\]
Its left adjoint is
\[
 M_0:\int_{\mathrm{Sm}_S}\SH(-)\longrightarrow\SH(S),
 \qquad (X,E)\longmapsto p_\sharp E,
\]
where $p:X\to S$ is the structure morphism
\cite[Construction~2.3]{mcolimits}. Indeed, a map $(X,E)\to(S,B)$
necessarily lies over $p$, so its mapping space is
\[
 \Map_{\SH(X)}(E,p^*B)\simeq\Map_{\SH(S)}(p_\sharp E,B).
\]

By adjunction, $M_0$ carries a canonical oplax symmetric monoidal
structure \cite[Corollary~7.3.2.7]{LurieHA}.
We check that its comparison maps are equivalences.
Let $q:Y\to S$ be the structure morphism and put
$r=p\circ\mathrm{pr}_X=q\circ\mathrm{pr}_Y$.
The binary comparison
\[
 r_\sharp(\mathrm{pr}_X^*E\otimes\mathrm{pr}_Y^*F)
 \longrightarrow p_\sharp E\otimes q_\sharp F
\]
Is the composite
By functoriality of smooth pushforward and
\cite[Theorem~4.10(i)(a)--(c)]{KR}, we have
\begin{align*}
& r_\sharp\bigl(\mathrm{pr}_X^*E\otimes\mathrm{pr}_Y^*F\bigr)\\
&\simeq p_\sharp(\mathrm{pr}_X)_\sharp
    \bigl(\mathrm{pr}_X^*E\otimes\mathrm{pr}_Y^*F\bigr)
    && \text{since } r=p\circ\mathrm{pr}_X\\
&\simeq p_\sharp\bigl(E\otimes
    (\mathrm{pr}_X)_\sharp\mathrm{pr}_Y^*F\bigr)
    && \text{by the projection formula for } \mathrm{pr}_X\\
&\simeq p_\sharp\bigl(E\otimes p^*q_\sharp F\bigr)
    && \text{by smooth base change for } q \text{ along } p\\
&\simeq p_\sharp E\otimes q_\sharp F
    && \text{by the projection formula for } p.
\end{align*}
The unit comparison is the canonical identification
$(\mathrm{id}_S)_\sharp\mathbf1_S\simeq\mathbf1_S$.
Thus $M_0$ is symmetric monoidal.

Now let $\int_{\mathrm{Sm}_S}\mathrm{Sph}$ denote the category of
elements of the presheaf $X\mapsto\Pic(\SH(X))$.
Its objects are pairs $(X,L)$ with $L$ tensor-invertible in $\SH(X)$.
A morphism $(X,L)\to(Y,N)$ over $h:X\to Y$ is an equivalence
$L\simeq h^*N$ in $\SH(X)$. The fiberwise Picard inclusions are symmetric monoidal and natural under
pullback. By symmetric monoidal unstraightening, they induce the symmetric
monoidal functor
\[
j:\int_{\mathrm{Sm}_S}\mathrm{Sph}
\longrightarrow\int_{\mathrm{Sm}_S}\SH(-).
\] 
Hence
$M_0\circ j$ is symmetric monoidal.

By \cite[Lemma~A.4(1)]{mcolimits}, the category of maps
$yX\to\mathrm{Sph}$ with $X\in\mathrm{Sm}_S$, denoted by
$(\mathrm{Sm}_S)_{/\mathrm{Sph}}$, is the right fibration over
$\mathrm{Sm}_S$ classified by $\mathrm{Sph}$. Thus unstraightening gives
\[
\int_{\mathrm{Sm}_S}\mathrm{Sph}
\simeq (\mathrm{Sm}_S)_{/\mathrm{Sph}}
\]
over $\mathrm{Sm}_S$, sending $(X,L)$ to the map $yX\to\mathrm{Sph}$
classified by $L$ under Yoneda.

Passing to presheaves and then applying \cite[Lemma~A.5]{mcolimits} gives
\[
\mathcal P\!\left(\int_{\mathrm{Sm}_S}\mathrm{Sph}\right)
\simeq \mathcal P\!\left((\mathrm{Sm}_S)_{/\mathrm{Sph}}\right)
\simeq \mathcal P(\mathrm{Sm}_S)_{/\mathrm{Sph}}.
\]
The represented families form a symmetric monoidal subcategory of
the operadic slice: Yoneda preserves finite products, and the
structure map of their tensor is the corresponding multiplication
in $\mathrm{Sph}$. This identifies their tensor with the external
tensor on $\int_{\mathrm{Sm}_S}\mathrm{Sph}$.

The slice tensor preserves colimits separately, since products
of presheaves of spaces preserve colimits separately and slice
colimits are computed on underlying presheaves. By
\cite[Proposition~4.8.1.10(4) and Corollary~4.8.1.12]{LurieHA},
the symmetric monoidal functor $M_0\circ j$ therefore extends to
a colimit-preserving symmetric monoidal functor on
$\mathcal P(\mathrm{Sm}_S)_{/\mathrm{Sph}}$. Its underlying functor
is $M$ by \cite[Definition~16.1 and Remark~16.5]{norms}.

\end{proof}
The motivic Thom functor preserves colimits and its source and target
are presentable. Write
\[
 M:\mathcal P(\mathrm{Sm}_S)_{/\mathrm{Sph}}
 \rightleftarrows \SH(S):M^R
\]
for the resulting adjunction. By
Lemma~\ref{lem:thom-families-monoidal}, $M$ is symmetric monoidal
for the tensor product of families. Hence $M^R$ is canonically lax
symmetric monoidal and the adjunction induces an adjunction on
commutative algebras
\cite[Proposition~7.3.2.5 and Corollary~7.3.2.7]{LurieHA}.

By \cite[Proposition~A.5.3]{HKZ}, applied with $G=1$ and the commutative
operad:
\[
\CAlg\bigl(\mathcal P(\mathrm{Sm}_S)_{/\mathrm{Sph}}\bigr)
\simeq
\mathrm{CMon}\bigl(\mathcal P(\mathrm{Sm}_S)\bigr)_{/\mathrm{Sph}}.
\]

\begin{lemma}[Multiplicative lifts]
\label{lem:thom-multiplicative-lifts}
Let $A\in\CAlg(\SH(S))$, and let $f:I\to\mathrm{Sph}$ be a family of invertible
motivic spectra.
There is a natural equivalence
\[
 \Map_{\CAlg(\SH(S))}(M(f),A)
 \simeq
 \Map_{\mathrm{CMon}(\mathcal P(\mathrm{Sm}_S))_{/\mathrm{Sph}}}
     \bigl((I,f),M^R(A)\bigr).
\]
\end{lemma}

\begin{proof}
This is the mapping-space equivalence for the adjunction on
commutative algebras induced by $M\dashv M^R$.
\end{proof}

\begin{lemma}[The right adjoint of the motivic Thom functor]
\label{lem:thom-right-adjoint}
Let $A\in\CAlg(\SH(S))$, and let $f:I\to\mathrm{Sph}$ be a family of invertible
motivic spectra. Under straightening, the pullback

$f^*M^R(A)\to I$ corresponds to the presheaf
\[
 \left(\int_{\mathrm{Sm}_S}I\right)^{\mathrm{op}}
 \longrightarrow \mathcal S,
 \qquad
 (X,i)\longmapsto\Map_{\SH(X)}\bigl(f_X(i),A_X\bigr).
\]
In particular, the fiber of $M^R(A)(X)\to\mathrm{Sph}(X)$ at $L$
is naturally equivalent to $\Map_{\SH(X)}(L,A_X)$.

\end{lemma}

\begin{proof}
Let $yX$ be the presheaf represented by $p:X\to S$, and let
$\underline L:yX\to\mathrm{Sph}$ classify $L\in\mathrm{Sph}(X)$.
There are natural equivalences
\[
\begin{aligned}
 \operatorname{fib}_L\bigl(M^R(A)(X)\to\mathrm{Sph}(X)\bigr)
 &\simeq
 \Map_{\mathcal P(\mathrm{Sm}_S)_{/\mathrm{Sph}}}
      \bigl((yX,\underline L),M^R(A)\bigr)\\
 &\simeq \Map_{\SH(S)}\bigl(M(\underline L),A\bigr)\\
 &\simeq \Map_{\SH(S)}(p_\sharp L,A)\\
 &\simeq \Map_{\SH(X)}(L,A_X).
\end{aligned}
\]
The first equivalence is Yoneda in the slice, the second is the
adjunction $M\dashv M^R$, the third is the represented-family
formula \cite[Definition~16.1 and Remark~16.5]{norms}, and the
last is the adjunction $p_\sharp\dashv p^*$.

A morphism $(X',L')\to(X,L)$ in
$\int_{\mathrm{Sm}_S}\mathrm{Sph}$ consists of an $S$-morphism
$h:X'\to X$ and an equivalence $u:L'\simeq h^*L$. Under the
displayed identification, restriction sends $\alpha:L\to A_X$
to the composite
\[
 L'\xrightarrow{u}h^*L\xrightarrow{h^*\alpha}
 h^*A_X\simeq A_{X'}.
\]
Thus $M^R(A)$ corresponds under straightening to the presheaf
$(X,L)\mapsto\Map_{\SH(X)}(L,A_X)$
\cite[Proposition~2.9 and Lemma~A.5]{mcolimits}.
Pullback along $f$ corresponds to restriction along
\[
 \int_{\mathrm{Sm}_S}I\longrightarrow
 \int_{\mathrm{Sm}_S}\mathrm{Sph},
 \qquad (X,i)\longmapsto(X,f_X(i)),
\]
which gives the asserted formula for $f^*M^R(A)$. Its sections
are therefore the coherent families of maps $f_X(i)\to A_X$,
that is, the transformations $f\to\pi^*A$ for $\pi:I\to *$.
\end{proof}

\begin{lemma}[Multiplication on the right adjoint]
\label{lem:thom-right-adjoint-product}
Let $A\in\CAlg(\SH(S))$ and let $p:X\to S$ be in $\mathrm{Sm}_S$.
The product of a finite family of points
$(L_j,\alpha_j:L_j\to A_X)_{j\in J}$ of $M^R(A)(X)$ is

\[
 \left(
   \bigotimes_{j\in J}L_j,\quad
   \bigotimes_{j\in J}L_j
   \xrightarrow{\otimes_j\alpha_j}
   A_X^{\otimes J}\xrightarrow{\mu_{A_X}^{(J)}}A_X
 \right).
\]
Equivalently, under the free-module adjunction this product is the
relative tensor product of the corresponding $A_X$-module maps
$A_X\otimes L_j\to A_X$. The unit corresponds to
$\mathrm{id}_{A_X}$.

\end{lemma}

\begin{proof}
Let $A\in\CAlg(\SH(S))$. For each $j\in J$, let
$\bar\alpha_j:p_\sharp L_j\to A$ be adjoint to $\alpha_j$.
By definition of the lax symmetric monoidal structure on $M^R$,
the external product of the given points is adjoint to
\[
 M\!\left(\bigotimes_{j\in J}\underline L_j\right)
 \simeq \bigotimes_{j\in J}p_\sharp L_j
 \xrightarrow{\otimes_j\bar\alpha_j}
 A^{\otimes J}\xrightarrow{\mu_A^{(J)}}A.
\]
By the tensor comparison in
Lemma~\ref{lem:thom-families-monoidal}, the corresponding map
over $X_S^J$ is the tensor product of the maps
$\mathrm{pr}_j^*\alpha_j$ followed by multiplication in
$A_{X_S^J}$. Pulling back along the diagonal $X\to X_S^J$
gives the stated product in $M^R(A)(X)$. For $J=\varnothing$,
the same formula gives $(\mathbf1_X,\eta_{A_X})$.

Finally, the free-module adjunction identifies $\alpha:L\to A_X$
with the $A_X$-module map
\[
 \widetilde\alpha:
 A_X\otimes L\xrightarrow{1\otimes\alpha}
 A_X\otimes A_X\xrightarrow{\mu_{A_X}}A_X.
\]
Since the free-module functor is symmetric monoidal, the adjunct of
the $J$-ary product is
\[
 A_X\otimes\bigotimes_{j\in J}L_j
 \simeq \bigotimes_{j\in J}^{A_X}(A_X\otimes L_j)
 \xrightarrow{\otimes_{j\in J}^{A_X}\widetilde\alpha_j}
 \bigotimes_{j\in J}^{A_X}A_X
 \simeq A_X,
\]
where $\bigotimes^{A_X}$ denotes the tensor product in
$\Mod_{A_X}(\SH(X))$. In arity zero this is $\mathrm{id}_{A_X}$.
More precisely, write $\Mod_{A_X}$ for $\Mod_{A_X}(\SH(X))$.
The assignment $\alpha\mapsto\widetilde\alpha$ is induced by the
following composite of symmetric monoidal functors:
\[
 \SH(X)_{/A_X}
 \xrightarrow{\,A_X\otimes-\,}
 (\Mod_{A_X})_{/(A_X\otimes A_X)}
 \xrightarrow{\,\mu_{A_X}\circ-\,}
 (\Mod_{A_X})_{/A_X}.
\]
The first functor is induced by symmetric monoidal extension of
scalars. Here $A_X$ acts on the first factor of $A_X\otimes A_X$.
The second functor is postcomposition with the multiplication
$A_X\otimes A_X\to A_X$, the counit at the unit module of the
induced adjunction on commutative algebras
\cite[\S2.2.2, Corollary~4.2.4.8, and Remark~4.5.3.2]{LurieHA}.
These constructions are natural under every $S$-morphism $h:X'\to X$:
pullback carries the unit, multiplication, and module actions to
their counterparts over $X'$, and preserves relative tensor products
because it preserves the two-sided bar construction and its geometric
realization \cite[Construction~4.4.2.7 and Theorem~4.4.2.8]{LurieHA}.

\end{proof}
For $A\in\CAlg(\SH(S))$ we will denote
\[
 \Pic(\Mod_A)(X):=\Pic\bigl(\Mod_{A_X}(\SH(X))\bigr),
\]
where the tensor product of modules is $\otimes_{A_X}$.

Let $(M^R(A))^\times$ denote the presheaf of units of $M^R(A)$,
obtained at each $X$ by taking the union of the components
representing units in $\pi_0(M^R(A)(X))$. 
\begin{lemma}[Units are module trivializations]
\label{lem:thom-units-trivializations}
There is a natural
equivalence of presheaves of grouplike $\mathbb E_\infty$-spaces
over $\mathrm{Sph}$
\[
 (M^R(A))^\times
 \simeq
 \mathrm{Sph}\times_{\Pic(\Mod_A)}*,
\]
where $\mathrm{Sph}\to\Pic(\Mod_A)$ sends $L$ to $A_X\otimes L$
and $*\to\Pic(\Mod_A)$ selects the unit module $A_X$.
\end{lemma}

\begin{proof}
Fix $X\in\mathrm{Sm}_S$. By Lemmas~\ref{lem:thom-right-adjoint}
and~\ref{lem:thom-right-adjoint-product}, a point of $M^R(A)(X)$ consists of
an invertible spectrum $L\in\SH(X)$ and
an $A_X$-module map $\tau:A_X\otimes L\to A_X$. Multiplication
tensors these maps over $A_X$, and the unit corresponds to
$\mathrm{id}_{A_X}$. Thus forgetting $L$ defines a symmetric
monoidal functor
\[
 M^R(A)(X)\longrightarrow
 \bigl(\Mod_{A_X}(\SH(X))\bigr)_{/A_X},
 \qquad (L,\tau)\longmapsto \tau.
\]

Suppose $(L,\tau)$ is a unit. Then $\tau$ is tensor-invertible
in the slice. The tensor unit of the slice is $\mathrm{id}_{A_X}$,
which is also its terminal object. Tensoring by $\tau$ is an
equivalence of the slice and carries $\mathrm{id}_{A_X}$ to $\tau$.
Hence $\tau$ is terminal. The canonical map from $\tau$ to
$\mathrm{id}_{A_X}$ is therefore an equivalence in the slice. On underlying modules,
this map is precisely $\tau:A_X\otimes L\to A_X$.

Conversely, the subspace of $M^R(A)(X)$ consisting of pairs for
which $\tau$ is an equivalence is the homotopy pullback
\[
 \mathrm{Sph}(X)\times_{\Pic(\Mod_A)(X)}\{A_X\}.
\]
Indeed, a point of this pullback consists of $L$ together with a
chosen equivalence $A_X\otimes L\simeq A_X$. Its inclusion in
$M^R(A)(X)$ respects multiplication, since both products tensor
the module maps over $A_X$. The pullback is grouplike, so every
such pair is a unit of $M^R(A)(X)$.

Finally, restriction along any $S$-morphism $h:X'\to X$ sends
$(L,\tau)$ to $(h^*L,h^*\tau)$, using the canonical identification
$h^*(A_X\otimes L)\simeq A_{X'}\otimes h^*L$. Hence the inclusion
of the pullback defines a map of presheaves of commutative monoids
over $\mathrm{Sph}$. We have shown that it identifies its source
with the units at every $X$, giving the asserted equivalence.
\end{proof}

\subsection{Proof of Proposition \ref{prop:orientation}}\label{proofof}
\begin{proof}[Proof of Proposition~\ref{prop:orientation}]
By Lemma~\ref{lem:thom-multiplicative-lifts}, a commutative algebra
map $M(f)\to A$ is equivalent to a lift of $f:I\to\mathrm{Sph}$
through $M^R(A)\to\mathrm{Sph}$ in presheaves of commutative monoids.
Taking units is right adjoint to the inclusion of grouplike
commutative monoids. Since $I$ and $\mathrm{Sph}$ are grouplike,
applying this adjunction objectwise and taking the homotopy fiber
over $f$ identifies the spaces of lifts. Hence

\[
\begin{aligned}
 \Map_{\CAlg(\SH(S))}(M(f),A)
 &\simeq
 \Map_{\mathrm{CGrp}(\mathcal P(\mathrm{Sm}_S))_{/\mathrm{Sph}}}
      \bigl((I,f),(M^R(A))^\times\bigr)\\
 &\simeq
 \Map_{\mathrm{CGrp}(\mathcal P(\mathrm{Sm}_S))_{/\mathrm{Sph}}}
      \bigl((I,f),\mathrm{Sph}\times_{\Pic(\Mod_A)}*\bigr),
\end{aligned}
\]
where the second equivalence is
Lemma~\ref{lem:thom-units-trivializations}.

By the universal property of the pullback, the last mapping space
is the space of paths from
\[
 f_A:I\xrightarrow{f}\mathrm{Sph}\longrightarrow\Pic(\Mod_A)
 \qquad\text{to}\qquad
 I\longrightarrow *\xrightarrow{A}\Pic(\Mod_A)
\]
in $\Map_{\mathrm{CGrp}(\mathcal P(\mathrm{Sm}_S))}
       (I,\Pic(\Mod_A))$.
The second map selects the unit module $A_X$ at every
$X\in\mathrm{Sm}_S$. Thus these paths are precisely the null
homotopies of $f_A$ as a map of presheaves of grouplike
$\mathbb E_\infty$-spaces.
\end{proof}

For $\MGL_S$, the family in Proposition~\ref{prop:orientation} is the restriction of $J$ to virtual bundles of rank zero.

\begin{corollary}
An orientation $\MGL_S\to R$ is equivalent to a null homotopy of the composite
\[
 K^\circ\xrightarrow{e}K\xrightarrow{J}\mathrm{Sph}
 \longrightarrow\Pic(\Mod_R).
\]
\end{corollary}
\begin{proof}
Apply Proposition~\ref{prop:orientation} to $J\circ e$ and use the Thom identification and its multiplicative structure in \cite[Theorem~16.13 and proof of Theorem~16.19]{norms}.
\end{proof}

\section{\texorpdfstring{$\mathbb{G}_m$}{Gm}-Equivariant Motivic Spectra}
\label{sec:equivariant}

To define the motivic Tate construction, we work with
$\mathbb G_m$-equivariant motivic spectra, the motivic analogue of
$\mathbb T$-equivariant spectra. In this section and the next,
$G=\mathbb G_m$ acts trivially on the base scheme.

\subsection{Definitions}

\begin{definition}[$\mathbb G_m$-Spaces, {\cite[Definitions~3.3, 3.7 and~3.12]{equivariant}}]
The category $\mathcal H^G(S)$ of motivic $G$-spaces over $S$ is
the full subcategory of presheaves of spaces $F$ on smooth
$G$-quasi-projective $S$-schemes satisfying the following conditions.
\begin{itemize}
  \item (Homotopy invariance) Let $E\to X$ be a $G$-vector bundle
  and $p:Y\to X$ a $G$-equivariant torsor under $E$. Then $p^*:F(X)\to F(Y)$ is an equivalence.

  \item (Nisnevich excision) We have $F(\varnothing)\simeq *$.
  If $V\to X$ is equivariant and étale and induces an isomorphism
  $V\times_X Z\to Z$ over a $G$-invariant finitely presented closed
  subscheme $Z\subset X$, then the following square is cartesian:
\[
\begin{tikzcd}
F(X) \arrow[r] \arrow[d] & F(X\setminus Z) \arrow[d] \\
F(V) \arrow[r] & F(V\setminus (V\times_X Z)).
\end{tikzcd}
\]
\end{itemize}
\end{definition}

We will now define the equivariant sphere associated to a representation:

\begin{definition}[Representation Spheres, {\cite[\S5.2]{equivariant}}]
Let $V$ be a finite locally free $G$-representation over $S$, with
total space $\mathbb V(V)=\operatorname{Spec}_S\operatorname{Sym}(V^\vee)$
and zero section $0_S$. Its \emph{representation sphere} is the
pointed motivic $G$-space
\[
    S^V=\mathbb V(V)/(\mathbb V(V)\setminus 0_S).
\]
\end{definition}

The equivalence $S^{V\oplus W}\simeq S^V\wedge S^W$ identifies
direct sums of representations with smash products of their spheres.
We pass to the stable category by making suspension by each
representation sphere invertible.

\begin{definition}[$\mathbb G_m$-Spectra, {\cite[Definition~6.1 and Corollary~6.7]{equivariant}}]
The category $\mathcal{SH}^G(S)$ of genuine $G$-equivariant motivic
spectra over $S$ is obtained from the category $\mathcal H^G_*(S)$
of pointed motivic $G$-spaces by formally inverting the representation
spheres under smash product.
\end{definition}

The above definition has an interpretation in terms of motivic
spectra over the quotient stack $[S/G]$. More precisely, when
$S$ is $G$-quasi-projective, Khan and Ravi identify Hoyois's
equivariant motivic category with the motivic category over $[S/G]$
by a canonical symmetric monoidal equivalence
\[
    \mathcal{SH}^G(S)\simeq\mathcal{SH}([S/G])
\]
\cite[Remark~4.9]{KR}.
Hence we will use the
equivariant and quotient-stack descriptions interchangeably.
Since $G$ acts trivially on $S$, we have $[S/G]\simeq S\times BG$.
Write
\[
    q\colon S\times BG\longrightarrow S
\]
for the projection. Pullback along $q$ equips a motivic spectrum
over $S$ with the trivial $G$-action. Its right adjoint $q_*$ takes
genuine categorical fixed points \cite[Remark~6.20]{equivariant}.
\subsection{The Motivic Tate Construction}

We will now define the euler map:

\begin{definition}[Euler Maps, {\cite[Construction~8.1]{KR}}]
Let $V$ be a finite locally free $G$-representation over $S$.
The zero section of $\mathbb V(V)$ followed by the quotient map
induces a map
\[
    e_V\colon \mathbf 1_S\longrightarrow S^V
\]
in $\mathcal{SH}^G(S)$, called the \emph{Euler map} of $V$.
\end{definition}
\begin{discussion}\label{euler+}
Let
$\mathrm{Vect}^G_{\mathrm{emb}}$ be the symmetric monoidal $\mathcal P(\mathrm{Sm}_S)$-category
whose value at $U$ is the category of $G$-vector bundles on $U$ and
equivariant subbundle embeddings, with direct sum as tensor product.
Restriction along $U'\to U$ is pullback of bundles.

The zero bundle is initial in every context and is preserved by
pullback. Hence the inclusions $0\hookrightarrow V$ define an internal
natural transformation $\underline{0}\to\mathrm{id}$ on
$\mathrm{Vect}^G_{\mathrm{emb}}$. This transformation is monoidal,
since the direct sum of the inclusions into $V$ and $W$ is the
inclusion into $V\oplus W$.

The representation sphere construction gives a symmetric monoidal
$\mathcal P(\mathrm{Sm}_S)$-functor
\[
    S^{(-)}\colon \mathrm{Vect}^G_{\mathrm{emb}}
    \longrightarrow \SH^G(-).
\]
Indeed, subbundle embeddings preserve the complements of the zero
sections, the defining quotients commute with pullback, and direct
sums give smash products. Applying this $\mathcal P(\mathrm{Sm}_S)$-functor to
$\underline{0}\to\mathrm{id}$ gives the \emph{Euler transformation}
\[
    e_{(-)}\colon \underline{\mathbf1}\longrightarrow S^{(-)}.
\]
Its component at a $G$-vector bundle $V$ on $U$ is the Euler map
$e_V$. This is the motivic analogue of the Euler transformation of
Carmeli and Luecke~\cite[Definition~3.11]{SL}.
\end{discussion}
Let $\mathcal L=\mathcal O(1)$ be the tautological line bundle on
$\mathbb P_S^\infty$. For an oriented motivic ring spectrum $R$,
we would like to invert its first Chern class
$t=c_1^R(\mathcal L)\in R^{2,1}(\mathbb P_S^\infty)$.
To do this at the level of ring spectra, we need a ring spectrum
over $S$ whose cohomology records the $R$-cohomology of
$\mathbb P_S^\infty$. Borel cochains provide this ring:

\begin{definition}[Borel completion and Borel cochains]
Let $D\in\SH^G(S)$. Its \emph{Borel completion} $D\to\widehat D$
is the localization at maps that become equivalences after
forgetting the $G$-action. Its \emph{Borel cochains} are
\[
    D^{hG}:=q_*\widehat D.
\]
\end{definition}
We will now define the smashing localization inverting $e_{\mathcal L}$  via a telescope:

\begin{definition}[Inverting the Euler map, {\cite[Lemma~12.1]{norms}}]
For $E\in\SH^G(S)$, define
\[
    E[e_{\mathcal L}^{-1}]
    :=\operatorname*{colim}\bigl(
       E\longrightarrow E\otimes S^{\mathcal L}
        \longrightarrow E\otimes S^{2\mathcal L}
        \longrightarrow\cdots\bigr).
\]
Here $m\mathcal L$ denotes the direct sum of $m$ copies of
$\mathcal L$, and each transition map is given by tensoring
with $e_{\mathcal L}$.
\end{definition}
Since tensor product preserves colimits,
$E[e_{\mathcal L}^{-1}]\simeq E\otimes\mathbf1[e_{\mathcal L}^{-1}]$.

We can now finally define the Tate construction:
\begin{definition}[Motivic Tate construction]
\label{def:motivic-tate}
The \emph{Tate construction} of $D\in\SH^G(S)$ is
\[
    D^{tG}:=q_*\bigl(\widehat D[e_{\mathcal L}^{-1}]\bigr).
\]
\end{definition}
\begin{remark}
For a commutative algebra $D$, the spectrum $D^{tG}$ is again a
commutative algebra, since Borel completion and $q_*$ are lax
symmetric monoidal and the telescope is a smashing localization
\end{remark}

We give now a useful description of the borel completion:
\begin{lemma}
\label{lem:cochains}
Let $s_n:\mathbb P_S^n\to[S/G]$ classify $\mathcal O(1)$. For an equivariant motivic spectrum $D$ and an invertible equivariant spectrum $W$, there are natural equivalences
\[
 \widehat D\simeq\lim_n(s_n)_*s_n^*D,
 \qquad
 \widehat D\otimes W\simeq\lim_n(s_n)_*s_n^*(D\otimes W).
\]
For a commutative algebra $D$, the second equivalence is compatible with tensor products of invertible objects, using the relative tensor product over $\widehat D$.
\end{lemma}
\begin{proof}
Write $u:S\to[S/G]$ for the atlas and $p_n:\mathbb P_S^n\to S$
for the projection. The adjunction units give a map
\[
    D\longrightarrow \lim_n (s_n)_*s_n^*D.
\]
We show that its target is Borel complete and that pulling it back
along $u$ gives an equivalence.

We will now prove that $\lim_n (s_n)_*s_n^*D$ is borel complete:
Let $f$ be a map such that $u^*f$ is an equivalence. On each
standard chart $U\subset\mathbb P_S^n$, the trivialization of
$\mathcal O(1)$ identifies $s_n|_U$ with $u\circ(p_n|_U)$. Hence
\[
    (s_n^*f)|_U\simeq (p_n|_U)^*u^*f.
\]
Thus $s_n^*f$ is an equivalence by Zariski descent.
Hence \[
\begin{aligned}
\operatorname{Map}_{\mathcal{SH}^G(S)}
   \bigl(B,(s_n)_*s_n^*D\bigr)
&\overset{\mathrm{adj}}{\simeq}
  \operatorname{Map}_{\mathcal{SH}(\mathbb P_S^n)}
   (s_n^*B,s_n^*D)\\
&\xrightarrow[\simeq]{\,g\mapsto g\circ s_n^*f\,}
  \operatorname{Map}_{\mathcal{SH}(\mathbb P_S^n)}
   (s_n^*A,s_n^*D)\\
&\overset{\mathrm{adj}}{\simeq}
  \operatorname{Map}_{\mathcal{SH}^G(S)}
   \bigl(A,(s_n)_*s_n^*D\bigr).
\end{aligned}
\] 
Therefore each $(s_n)_*s_n^*D$ is Borel
complete, and so is their limit.

The pullback of $s_n$ along $u$ is the projection
$r_n:\mathbb A_S^{n+1}\setminus 0_S\to S$. Since $u$ is smooth,
$u^*$ preserves limits, and smooth base change gives
\[
    u^*\Bigl(\lim_n(s_n)_*s_n^*D\Bigr)
    \simeq \lim_n(r_n)_*r_n^*u^*D
    \qquad\text{\cite[Theorem~4.10(i)]{KR}}.
\]
The right-hand side is the spectrum of $u^*D$-valued cochains on
\[
    EG_S=\operatorname*{colim}_n
    (\mathbb A_S^{n+1}\setminus 0_S)
\]
over $S$. After forgetting the action, $EG_S$ is
$\mathbb A^1$-contractible over $S$: linear interpolation gives
homotopies from the identity to the coordinate shift and from the
shift to the constant first basis vector. The map from $u^*D$ to
these cochains is therefore an equivalence. This identifies the
original map with the Borel completion of $D$.

Finally, tensoring with $W$ preserves limits since $W$ is invertible.
The projection formula identifies each
$(s_n)_*s_n^*D\otimes W$ with $(s_n)_*s_n^*(D\otimes W)$, giving
the second equivalence. For commutative $D$, these identifications
respect multiplication, so the tensor comparison for $W$ and $W'$
is the canonical equivalence
\[
    (\widehat D\otimes W)\otimes_{\widehat D}(\widehat D\otimes W')
    \simeq \widehat D\otimes W\otimes W'.
\]
\end{proof}

For $R$ over $S$, write $\widehat R$ for the Borel completion of
$q^*R$ and $R^{hG}=q_*\widehat R$.
Applying $q_*$ to the first equivalence of Lemma~\ref{lem:cochains}
identifies $R^{hG}$ with the limit of the $R$-valued cochains
on $\mathbb P_S^n$. Thus $R^{hG}$ computes cochains on
$\mathbb P_S^\infty$, the motivic classifying space of line bundles
over $S$. In particular, the first Chern class of the tautological line bundle
is a class
\[
    t=c_1^R(\mathcal L)\in R^{2,1}(\mathbb P_S^\infty)
    \cong [\mathbf1_S,\Sigma^{2,1}R^{hG}]_{\SH(S)}.
\]
We now use the orientation to identify the Euler telescope, after
taking fixed points, with the telescope of multiplication by $t$.
\begin{proposition}
\label{prop:tate}
Let $R$ be a motivic $\mathbb E_\infty$-ring over $S$ with an
$\mathbb E_\infty$-orientation $\omega$. For each smooth $S$-scheme
$X$.
There are natural
equivalences of motivic $\mathbb E_\infty$-rings over $X$
\[
    (R^{tG})_X
    \simeq q_*\bigl(\widehat R_X[e_{\mathcal L}^{-1}]\bigr)
    \simeq R_X^{hG}[t^{-1}].
\]
\end{proposition}
\begin{proof}
We first take $X=S$.
Applying the lax symmetric monoidal functor $q_*$ to Euler
localization gives a map of commutative algebras
\[
    R^{hG}=q_*\widehat R
    \longrightarrow q_*\bigl(\widehat R[e_{\mathcal L}^{-1}]\bigr).
\]
Since $q_*$ preserves filtered colimits \cite[Theorem~4.10(iii)]{KR},
the target is computed by applying $q_*$ to each stage of the Euler
telescope. On each $\mathbb P_S^n$, the orientation $\omega$ gives
an $R$-linear Thom equivalence for $m\mathcal O(1)$, compatible with
restriction in $n$ and direct sums. Lemma~\ref{lem:cochains}
therefore gives an equivalence of $R^{hG}$-modules
\[
    q_*(\widehat R\otimes S^{m\mathcal L})
    \simeq \Sigma^{2m,m}R^{hG}.
\]
The transitions are the successive shifts of multiplication by $t$,
since the zero section pulls the Thom class of $\mathcal O(1)$
back to $c_1^R(\mathcal O(1))$.

Thus, on underlying modules, the map from $R^{hG}$ to the
telescope is localization at $t$
\cite[Lemma~12.1 and Example~12.4]{norms}. In particular, $t$ acts
invertibly on the target, so the universal property of localization
gives a map of commutative $R^{hG}$-algebras
\[
    R^{hG}[t^{-1}]
    \longrightarrow q_*\bigl(\widehat R[e_{\mathcal L}^{-1}]\bigr).
\]
Its underlying module map is an equivalence by uniqueness of
localization, so it is an equivalence of commutative algebras.

Finally, let $p:X\to S$ be smooth. Smooth pullback preserves the Borel
limit and the Euler telescope. Smooth base change for the direct images
\cite[Theorem~4.10(i)]{KR} identifies $p^*(R^{tG})$ with $(R_X)^{tG}$
and $p^*(R^{hG})$ with $R_X^{hG}$. Since $p^*t$ is the first Chern
class of the pulled-back tautological line bundle, this gives the stated
equivalences over $X$, naturally in $X$.
\end{proof}
We now define the equivariant motivic $J$-homomorphism by giving
each vector bundle the scalar $G$-action and taking its Thom sphere.

\begin{construction}[Equivariant motivic $J$-homomorphism]
Let $V$ be a vector bundle on $X$. With $G=\mathbb G_m$ acting by
scalar multiplication on its fibers, $V$ corresponds to the bundle
$q^*V\otimes\mathcal L$ on $[X/G]$. This gives the composite
\[
    \widetilde J_X\colon K(X)
    \xrightarrow{\;[\mathcal L]\cdot q^*(-)\;}
    K([X/G])
    \xrightarrow{\;J\;}
    \Pic(\SH^G(X)),
\]
which sends $[V]$ to $S^{q^*V\otimes\mathcal L}$.

\end{construction}

\section{Motivic Euler--Tate Characteristic Classes}
\label{sec:classes}

In this section we define the motivic Euler--Tate characteristic
class. We first introduce projective orientations and show that an
oriented motivic ring spectrum $(R,\omega)$ gives rise to two of
them. The Euler--Tate characteristic class is the quotient of the
projective orientation induced by $\omega$ by the Euler orientation.

\subsection{Projective Orientations in the Motivic Setting}

Recall that an orientation of $R$ gives, for each rank-zero virtual
bundle $\xi$ on $X$, a trivialization
\[
    R_X\otimes S^\xi \simeq R_X.
\]
We now give the bundle the scalar action. Its Thom sphere is then
$S^{q^*\xi\otimes\mathcal L}$, an equivariant motivic spectrum over $X$.

A projective orientation asks for a trivialization of this equivariant
Thom sphere after a chosen change of coefficients. In our case, we
will change coefficient to $\widehat R_X[e_{\mathcal L}^{-1}]$. The trivialization we seek is
therefore an equivalence of modules
\[
    \widehat R_X[e_{\mathcal L}^{-1}]
    \otimes S^{q^*\xi\otimes\mathcal L}
    \simeq \widehat R_X[e_{\mathcal L}^{-1}].
\]
In the general definition below, the symmetric monoidal functor
$I_X$ specifies the change of coefficients.

\begin{definition}[Projective orientations]
Let $I\colon\SH^G(-)\to\mathcal C(-)$ be a symmetric monoidal
functor parametrized over $\mathrm{Sm}_S$, and write $I_X$ for
its component at $X$. An
\emph{$(\mathbb E_\infty,I)$-projective orientation} is a null
homotopy of the map of presheaves of grouplike
$\mathbb E_\infty$-spaces given over $X$ by
\[
 K^\circ(X)\longrightarrow K(X)\xrightarrow{\ \widetilde J_X\ }
 \Pic(\SH^G(X))\xrightarrow{\ \Pic(I_X)\ }\Pic(\mathcal C(X)).
\]
Here we take $I_X$ to be the free $\widehat R_X[e_{\mathcal L}^{-1}]$-module functor.
\end{definition}

Over each projective space in the Borel construction, the scalar
twist becomes tensoring with $\mathcal O(1)$. For a rank-zero
virtual bundle, the orientation $\omega$ of $R$ trivializes the
resulting Thom module. We now use these trivializations to construct
the first of our two projective orientations, $\omega_{\mathrm{prj}}$.
\begin{construction}
Let $R$ have an orientation $\omega$ and let $\xi\in K^\circ(X)$. For a smooth $X$-scheme $p:Y\to X$ equipped with a line bundle $L$, the class $p^*\xi\otimes L$ has rank zero. By Proposition~\ref{prop:orientation}, $\omega$ gives an equivalence
\[
 R_Y\otimes S^{p^*\xi\otimes L}\simeq R_Y
\]
natural in $(Y,L)$ and symmetric monoidal in $\xi$. Applied to the tautological line bundles $\mathcal O(1)$ on the projective spaces over $X$, these equivalences are compatible under restriction. By Lemma~\ref{lem:cochains} they give a symmetric monoidal trivialization of $\widehat R_X\otimes\widetilde J_X(\xi)$. Extension of scalars to $\widehat R_X[e_{\mathcal L}^{-1}]$ gives the projective orientation $\omega_{\mathrm{prj}}$.
This is the motivic analogue of the projective orientation that
Carmeli and Luecke construct from $\omega$ and then push forward
to the Tate quotient~\cite[Definition~3.10 and the beginning of \S4.1]{SL}.
\end{construction}

For the second trivialization we use the euler map itself.

\begin{construction}[Euler orientation]
Let $V$ be a vector bundle on $X$. The euler map of
$q^*V\otimes\mathcal L$ induces a map of
$\widehat R_X[e_{\mathcal L}^{-1}]$-modules
\[
    a_V\colon \widehat R_X[e_{\mathcal L}^{-1}]
    \longrightarrow
    \widehat R_X[e_{\mathcal L}^{-1}]
    \otimes S^{q^*V\otimes\mathcal L}.
\]
On an open $U$ where $V$ is trivial of rank $r$, the bundle
$q^*V\otimes\mathcal L$ is a direct sum of $r$ copies of
$\mathcal L$. Its Euler map is therefore $e_{\mathcal L}^{\otimes r}$.
Since we have inverted $e_{\mathcal L}$, the map $a_V$ is an
equivalence over $U$, and hence over $X$ by Zariski descent.

And the inverses $a_V^{-1}$ trivialize the Thom modules. From \ref{euler+} these
trivializations are natural under pullback and symmetric monoidal
under direct sum.
Since $\Pic(\mathcal C)$ is
grouplike, they extend to the group completion of the groupoid of
vector bundles, whose Zariski sheafification is $K$
\cite[\S16.2, proof of Lemma~16.12]{norms}.
The target $\Pic(\mathcal C)$ satisfies Zariski descent, thus the trivialization extends to $K$. Its restriction to $K^\circ$
is the Euler orientation $\omega_{\mathrm{eul}}$.
\end{construction}
Hence after extension of scalars to $\widehat R_X[e_{\mathcal L}^{-1}]$, we get the trivialization supplied by $\omega$ and the one supplied by the zero section.
Therefore each oriented rings is equipped with two orientations:
\begin{corollary}
Every oriented ring $R$ admits two natural projective orientations in the above Tate construction, $\omega_{\mathrm{prj}}$ and $\omega_{\mathrm{eul}}$.
\end{corollary}

\subsection{Main Definition and Results}

A projective orientation is a null homotopy of the fixed map
$K^\circ\to\Pic(\mathcal C)$. Thus, as in~\cite[Remark~3.8]{SL},
the space of projective orientations is a torsor under the
underlying space of the connective mapping spectrum
\[
\begin{aligned}
    \operatorname{Map}\bigl(K^\circ,\Omega_{\mathbf1}\Pic(\mathcal C)\bigr)
    &\simeq \operatorname{Map}\bigl(K^\circ,\mathrm{GL}_1(R^{t\mathbb G_m})\bigr)\\
    &\simeq \operatorname{Map}\bigl(\mathrm{BGL},\mathrm{GL}_1(R^{t\mathbb G_m})\bigr).
\end{aligned}
\]
For the first equivalence, loops at the unit are its automorphisms
\cite[\S2.2, before Example~2.2.2]{MathewStojanoska}, evaluation at $1$ identifies automorphisms of the unit module
$\widehat R_X[e_{\mathcal L}^{-1}]$ with its invertible sections
over $[X/G]$. Since $q^*\mathbf1_X\simeq\mathbf1_{[X/G]}$, the
adjunction $q^*\dashv q_*$ identifies these with of $q_*(\widehat R_X[e_{\mathcal L}^{-1}])$ over $X$.
By Proposition~\ref{prop:tate}, this ring is equivalent to
$(R^{t\mathbb G_m})_X$, so its space of units is
$\mathrm{GL}_1(R^{t\mathbb G_m})(X)$.
The second equivalence uses $\mathrm{BGL}=L_{\mathrm{mot}}K^\circ$
\cite[\S16.2, before Theorem~16.13]{norms} and that the presheaf of units
is motivically local
\cite[proof of Proposition~2.2.3]{MathewStojanoska}.
\\
Interchanging the two paths inverts the resulting unit, so we fix their order:

\begin{definition}[Quotient of Projective Orientations]
Given two projective orientations $\omega_1,\omega_2$, their
\emph{quotient} $\omega_1/\omega_2$ is the loop obtained by following
$\omega_2$ backwards and then following $\omega_1$. Thus
$(\omega_1/\omega_2)\circ\omega_2\simeq\omega_1$.
\end{definition}

We now define the motivic Euler--Tate characteristic class by dividing
the projective orientation induced by $\omega$ by the Euler orientation.

\begin{definition}[Euler--Tate Characteristic Class]
Let $R$ be a motivic $\mathbb E_\infty$-ring with an
$\mathbb E_\infty$-orientation $\omega$. Its \emph{Euler--Tate
characteristic class} is the map
\[
    c^{\mathrm{et}}_{R,\omega}
    :=\frac{\omega_{\mathrm{prj}}}{\omega_{\mathrm{eul}}}
    \colon \mathrm{BGL}\longrightarrow
          \mathrm{GL}_1(R^{t\mathbb G_m}).
\]
This is a map of grouplike motivic $\mathbb E_\infty$-spaces. For a rank-$r$ vector bundle $V$
on $X$, we write $c^{\mathrm{et}}_{R,\omega}(V)$ for its value on
the rank-zero class $[V]-[\mathcal O_X^r]$.
\end{definition}

We now compute the Euler--Tate characteristic class of a vector
bundle in terms of its Chern roots and the formal group law of the
chosen orientation.

\begin{theorem}[Evaluation on vector bundles]
\label{thm:evaluation}
Let $(R,\omega)$ be an $\mathbb E_\infty$-oriented motivic ring
spectrum, and let $V$ be a vector bundle of rank $r$ on
$X\in\mathrm{Sm}_S$.
Then
\[
    c^{\mathrm{et}}_{R,\omega}(V)
    =\frac{e_R(q^*V\otimes\mathcal L)}{t^r}
    \quad\text{in }(R^{t\mathbb G_m})^{0,0}(X)^\times.
\]
\end{theorem}

\begin{proof}
We may take $X=S$ by applying the calculation over $X$ with
coefficients $R_X$ and the pullback of $\omega$.
By definition, we evaluate the quotient
$\omega_{\mathrm{prj}}/\omega_{\mathrm{eul}}$ on
$[V]-[\mathcal O_S^{\oplus r}]$. Following the Euler orientation
backwards and then applying the orientation induced by $\omega$
amounts to composing the zero section with the Thom class for $V$ this give $e_R(q^*V\otimes\mathcal L)$.
And for 
$\mathcal O_S^{\oplus r}$ we get $e_R(\mathcal L)^r=t^r$. Since we evaluate on the virtual
difference, we divide the first contribution by the second.
\end{proof}
We recall that the \emph{characteristic series} of
$c^{\mathrm{et}}_{R,\omega}$ is the power series $f(x)$ with
coefficients in $(R^{t\mathbb G_m})^{*,*}(S)$ such that, for every
line bundle $L$ on a smooth $S$-scheme $X$,
\[
    f\bigl(c_1^R(L)\bigr)=c^{\mathrm{et}}_{R,\omega}(L)
    \quad\text{in }\bigl((R^{t\mathbb G_m})^{0,0}(X)\bigr)^\times.
\]
The coordinate $x=c_1^R(\mathcal O(1))$ is the first Chern class
of the universal line on $\mathbb P_S^\infty$. Its pullback is
$c_1^R(L)\in R^{2,1}(X)$, which we view in Tate cohomology via
$R\to R^{t\mathbb G_m}$ (cf.~\cite[\S2.2, Proposition~2.7]{SL}).
It is obtained by evaluating $c^{\mathrm{et}}_{R,\omega}$ on the
universal line bundle $\mathcal O(1)$ over $\mathbb P_S^\infty$
and writing the resulting Tate cohomology class as a power series
in $x=c_1^R(\mathcal O(1))$
(cf.~\cite[proof of Proposition~4.4]{SL}).

\begin{theorem}[Characteristic series, {cf.~\cite[Proposition~4.4]{SL}}]
\label{thm:characteristic-series}
Let $F_R$ be the formal group law of $\omega$. The characteristic
series of $c^{\mathrm{et}}_{R,\omega}$ is
\[
    \frac{F_R(x,t)}{t},
\]
where $x=c_1^R(\mathcal O(1))$ is the class of the universal line
on $\mathbb P_S^\infty$, viewed in Tate cohomology. Consequently,
on a splitting space for $V$ with Chern roots $x_1,\ldots,x_r$,
\[
    c^{\mathrm{et}}_{R,\omega}(V)
    =\prod_{j=1}^r\frac{F_R(x_j,t)}{t}.
\]
\end{theorem}

\begin{proof}
The characteristic series is the value of $c^{\mathrm{et}}_{R,\omega}$
on the universal line. Apply Theorem~\ref{thm:evaluation} to
$\mathcal O(1)$ on $\mathbb P_S^n$. The two lines in the numerator
come from different factors. The line $\mathcal O(1)$ on
$\mathbb P_S^n$ has first Chern class $x$, and the tautological line bundle
$\mathcal L$ from the Borel construction has first Chern class $t$.
Their tensor product therefore has first Chern class $F_R(x,t)$.
Dividing by $t$ and passing to the limit in $n$ gives the stated
series.

On the full flag bundle of $V$, the tautological filtration has
line quotients $L_1,\ldots,L_r$. In $K_0$ the pullback of
$[V]-[\mathcal O_X^{\oplus r}]$ is the sum of the classes
$[L_j]-[\mathcal O]$. Since $c^{\mathrm{et}}_{R,\omega}$ takes
sums to products, the formula follows with $x_j=c_1^R(L_j)$.
\end{proof}
For Spitzweck's integral motivic cohomology~\cite[Definition~4.27 and \S9]{Spitzweck} the formal group law is additive, so a line with first Chern class $x$ contributes the factor $1+xt^{-1}$.
\begin{corollary}
\label{cor:integral}
For Spitzweck's integral motivic cohomology spectrum $R=H\mathbb Z_S$ with its canonical orientation,
\[
    c^{\mathrm{et}}_{H\mathbb Z}(V)
      =\sum_{i=0}^{r} c_i(V)\,t^{-i}.
\]
\end{corollary}
\begin{proof}
The canonical map $\mathrm{MGL}_S\to H\mathbb Z_S$ is a map of
$\mathbb E_\infty$-rings \cite[Remark~4.3.3]{modules} and has additive
formal group law \cite[\S10, between Proposition~10.1 and Remark~10.2]{Spitzweck}.
The preceding product is therefore $\prod_j(1+x_jt^{-1})$.
Its coefficient of $t^{-i}$ is the $i$th elementary symmetric
polynomial in the Chern roots, hence $c_i(V)$.
\end{proof}

We now compute the bigraded coefficients of the motivic Tate construction.

\begin{proposition}
\label{prop:coefficients}
Let $R$ be an $\mathbb E_\infty$-oriented motivic ring spectrum. For all $a,b\in\mathbb Z$ and every smooth qcqs $X$ over $S$, the
bigraded groups of the Euler-localized Borel theory are
\[
 (R^{t\mathbb G_m})^{a,b}(X)
   \cong\underset{m\geq 0}{\operatorname{colim}}\;
      \prod_{j\geq -m} R^{a-2j,b-j}(X)\,t^j .
\]
The transition maps insert zero in the new lowest exponent, and
with the usual cohomological sign convention the multiplication is the
Cauchy product. In particular, writing
$H^{a,b}(S;\mathbb Z)=(H\mathbb Z_S)^{a,b}(S)$, one has
\[
 \pi_{p,q}\bigl((H\mathbb Z_S)^{t\mathbb G_m}\bigr)
 \cong\underset{m\geq 0}{\operatorname{colim}}\;
       \prod_{j\geq-m} H^{-p-2j,-q-j}(S;\mathbb Z)\,t^j .
\]
\end{proposition}
\begin{proof}
For the first formula it suffices to take $X=S$. Indeed, for general
$X$ we apply the calculation over $X$ to $R_X$ and use
$(R^{t\mathbb G_m})_X\simeq (R_X)^{t\mathbb G_m}$
from Proposition~\ref{prop:tate}.
The projective bundle formula for an oriented spectrum over a qcqs base
\cite[\S6, pp.~271--272]{AHI} gives
\[
 R^{a,b}(\mathbb P^n_S)
    =\bigoplus_{j=0}^{n} R^{a-2j,b-j}(S)\,t^j.
\]
For a noetherian base of finite Krull dimension this is also
\cite[Theorem~2.11]{NSO}. Restriction from $\mathbb P^{n+1}_S$ to
$\mathbb P^n_S$ is truncation since $t^{n+1}=0$ on $\mathbb P^n_S$.
Thus the inverse system is surjective in every bidegree. The Milnor
exact sequence \cite[Lemma~5.41]{ThomasonMilnor} identifies Borel cohomology in bidegree $(a,b)$ with the product over $j\geq 0$. Finally, the unit is compact in $\SH(S)$
\cite[Proposition~C.12(2)]{HoyoisTrace}, hence so are
its bidegree shifts. Thus
cohomology commutes with the telescope inverting $t$. Multiplication
by $t^{-m}$ shifts the lower bound to $-m$, giving the first formula.
For the second, use $\pi_{p,q}D=[S^{p,q},D]_{\SH(S)}$, so $(a,b)=(-p,-q)$.
\end{proof}

Over a field $k$, Spitzweck's spectrum agrees with Voevodsky's
\cite[proof of Theorem~B.4]{norms}, and motivic cohomology in
negative weights vanishes \cite[Corollary~2]{Chow}.
The corresponding bigraded coefficient ring in cohomological degrees
is therefore $H^{*,*}(k;\mathbb Z)[t,t^{-1}]$, with
$|t|=(2,1)$. In each bidegree the weight bounds the exponent above,
and Euler localization bounds it below.

Equivalently, for $d=p-2q$ one has
\[
 \pi_{p,q}\bigl((H\mathbb Z_k)^{t\mathbb G_m}\bigr)
  \cong \bigoplus_{s=0}^{d}\mathrm{CH}^{s}(k,d)\,t^{-q-s},
\]
with the empty sum interpreted as zero. Indeed,
$H^{2s-d,s}(k;\mathbb Z)=\mathrm{CH}^{s}(k,d)$ \cite[Corollary~2]{Chow}, and the
degree-$d$ simplex over $k$ has no cycles of codimension greater than $d$.
In particular, the diagonal $p=2q$ consists of $\mathbb Z\,t^{-q}$.

\end{document}